\documentclass[12pt]{article}
\usepackage[ansinew]{inputenc}
\usepackage[pdftex]{graphicx}
\usepackage{amssymb,amsmath,amscd}
\usepackage{dsfont}
\usepackage{verbatim}
\usepackage{amsmath} 
\usepackage{graphicx}
\usepackage{pstricks}
\usepackage{enumitem}
\usepackage{ntheorem}
\usepackage{amssymb}
\usepackage{float}
\usepackage{breqn}

\newtheorem{de}{Definition}[section] 
\newtheorem{nott}{Notation}

\newtheorem{theo}{Theorem}
  
\newtheorem{lem}{Lemma}[section]
\newtheorem{cor}[lem]{Corollary} 
\newtheorem{conj}{Conjecture}[section]

\setenumerate[1]{font=\bfseries,label=\Roman*.}
\setenumerate[2]{font=\itshape,label=\arabic*)}

\author{Yashar Memarian
}

\title{On the Maximum Number of Vertices of Critically Embedded Graphs \\
       }
\date{}

\begin{document}
\maketitle

\begin{abstract}
Define a boundary point of a graph which is embedded in the Euclidean plane a vertex which is incident to only one edge. In this paper we consider graphs which are embedded in the Euclidean plane with a finite number of boundary points. The simple geometric condition we impose on them is that the sum of unit vectors of edges extending from each non-boundary vertex will be equal to zero. We call such a graph a critical graph and ask to maximise the number of vertices of critical graphs with a given size of boundary. The main results of this paper give a sharp upper bound for the maximum number of vertices of planar critical graphs, where the degree of each non-boundary vertex is 3 or 4.
\end{abstract}

\section{Minimal and Critical Graphs}

Let $A$ be a finite number of points fixed on the Euclidean plane. We say that graph $G$ has the set $A$ as a boundary and denote it by $\delta G=A$, if $G$ contains {\em half-edges} (i.e. edges that are incident with only one end vertex) and the set of end vertices of all half-edges is $A$. The famous one-dimensional Plateau problem concerns the following variational problem: Consider all planar graphs $G$ such that $\delta G =A$ where $\delta G$ denotes the vertices (points) on the boundary of $G$ and define the functional $l:\{G\}\to\mathbb{R}_{+}$ by 
\begin{eqnarray*}
l(G)=\Sigma_{e\in E(G)}L(e),
\end{eqnarray*}
where $E(G)$ stands for the set of edges in $G$ and $L(e)$ denotes the (Euclidean) length of the edge $e$. Find the graph $G_{min}$ which minimises the functional $l$. Such a graph is called a minimal graph. It turns out that the edges of such a graph are straight lines (not very surprising) and the sum of unit vectors of edges extending from each non-boundary vertex is equal to zero. See \cite{Alm3} and references within for more discussions on the Plateau problem and see \cite{mark} for more details on minimal graphs. 

In this paper we are not interested in the length of graphs. Our graphs contain a boundary (i.e.the set of all end vertices of all the half-edges). The graphs are drawn in the plane so that each edge is a straight line. Edges are allowed to cross each other, and the intersecting point of two edges is a vertex of the graph. The (only) geometric condition we impose on our graph, is that the sum of unit vectors of edges extending from each non-boundary vertex must be equal to zero. We call these graphs \emph{critical}. Every minimal graph is critical-but the converse is not true.

The problem in which we are interested is the following : \emph{Find the maximum number of vertices of a critical graph with a given size of the boundary}. Note that there are considerable differences between our problem and the Plateau problem. In our problem, we do not consider the length of the edges at any point and we only require the \emph{size} of the boundary to be fixed. In the Plateau problem, the (actual) boundary points are considered to be fixed. Our question, simple as it appears, is completely open. If we do not impose some conditions on our graphs, this problem becomes completely pointless. For example, vertices with degree $2$ are not allowed, since we can add as many of these vertices to our graph and the resulting new graph remains critical. We also consider our graphs to be bounded, i.e. they are entirely contained in a bounded open set of the plane. And last, we impose that the graphs have countable vertices. From now on, the graphs studied in this paper will satisfy all of these characteristics. Keeping this in mind, it is unclear whether we can construct critical graphs with a finite number of boundary points and infinite number of vertices. In Section $4$ we shall formulate a few conjectures regarding this problem.

Here we study a rather simplified version of this question. In Section $2$ and $3$, we consider critical graphs for which the degree of each non-boundary vertex is  equal to $3$ (resp to $4$) and call them \emph{$3$-boundary regular critical graphs} (resp $4$-boundary regular critical graphs).

A few notations:



$G$ will denote the graph under study, $\vert G\vert$ denotes the number of vertices of $G$. $\delta G$ denotes the vertices of the boundary of $G$ which we refer to as boundary points all along this paper. The size of the boundary of $G$ denoted by $\vert \delta G\vert$ denotes the number of boundary vertices of $G$. For $n \geq 3$, we denote by $f_3(n)$ (resp. $f_4(n)$) the supremum of the
numbers of vertices of  $3$ (resp. $4$)-boundary regular critical graphs on the plane with $n$ boundary points.

Here are the two main theorems of this paper:

\begin{theo} \label{1}
Let $n$ be an integer and let $k=\lfloor n/6\rfloor$ and $l=n-6k$. The maximum number of vertices of a $3$-boundary regular critical graph in the plane with boundary of size $n$ is equal to 
\begin{eqnarray*}
f_3(n)=6k^2+6k+2l(k+1)-r,
\end{eqnarray*}
where $r=0$ if $n$ is a multiple of $6$ and $r=2$ otherwise. Furthermore, the combinatorial structure of critical graphs which maximise the number of vertices is unique.
\end{theo}


The combinatorial structure of the $3$-boundary regular critical graphs maximising the number of vertices with a given boundary size will be described in Definition \ref{defhn}.

\begin{theo} \label{sec}
Let $G$ be a $4$-boundary regular critical graph on the plane with boundary of size $n$. Then
\begin{eqnarray*}
f_4(n)=\binom{n/2}{2}+n  \qquad  \textrm{if n is even}\\
f_4(n)=0  \qquad \textrm{otherwise}.
\end{eqnarray*}
\end{theo}

Theorems \ref{1} and \ref{sec} are proven in Sections 2 and 3. Both proofs are elementary.

In the last section we will generalise this problem for the case of Riemannian manifolds and study a few examples. Finally, we will present a few open questions.

\section{$3$-Boundary Regular Critical Graphs}




\subsection{Preliminaries}


 

Let $G$ be a $3$-boundary regular critical graph with $n$ boundary points. By the definition of critical graphs, we know that the angle between the edges directing from every vertex is equal to $120$ degrees. This is actually the only geometric restriction put on the graph which will make the estimation of the function $f_3(n)$ easy.


\begin{lem} \label{increasing}
$f_3(n)$ is an increasing function of $n$.
\end{lem}

\emph{Proof of Lemma \ref{increasing}}:

Let $G$ be a $3$-boundary regular graph with $n$ boundary points. We choose one vertex $x \in \delta G$. Since all boundary points have a degree equal to $1$, there is only one edge $e$ directed from $x$ (i.e. the half edge $e$ with endvertex $x$). We add two edges directing from $x$ in a way that the angles between each of them and $e$ are equal to $120$ degrees. The new graph $G'$ is $3$- boundary regular and critical with $n+1$ boundary points, and $\vert G'\vert=\vert G\vert +2$. This completes the proof of the lemma.
\begin{flushright}
$\Box$
\end{flushright}

\begin{de}[Cycle, Interior]
A \emph{cycle} $C$ in $G$ is a $2$-regular subgraph of $G$. As a cycle is a simple closed curve in the plane, it separates the plane into two components. The bounded one is called the {\em interior} of $C$, and denoted by $\mathrm{Int}(C)$.
\end{de}

\begin{de}[Incoming and Outgoing Vertices and Edges] \label{in}
Let $C$ be a cycle in $G$. If $v\in V(C)$ is a vertex on $C$, let $e$ be the edge incident with $v$ that does not belong to $C$. If $e$ is embedded in $Int(C)$, then we say that $v$ and $e$ are \emph{incoming} vertices and edges respectively. Otherwise, they are \emph{outgoing}.
We denote by $V_{out}^{C}$ the number of outgoing vertices of the cycle $C$, and by $V_{in}^{C}$ the number of incoming vertices of the cycle $C$.
\end{de}

We define a partial order on the set of cycles.

\begin{de}[Maximal cycles] \label{defmax}
Let $C$ and $C'$ be two cycles. We define $C \leq C'$ if $\mathrm{Int}(C)\subseteq \mathrm{Int}(C')$. We call a cycle $C$ of $G$ \emph{maximal} if $C$ is maximal for this partial order.
\end{de}


\begin{lem} \label{out}
Let $C$ be a maximal cycle. For every outgoing vertex $v$ of $C$ the outgoing edge attached to $v$ does not belong to any cycle.
\end{lem}
\emph{Proof of Lemma \ref{out}}:

By contradiction, assume there is an outgoing edge $e$ of $C$ which belongs to a cycle $C'$. The cycles $C$ and $C'$ will have some edges in common. The outer boundary component of $\mathrm{Int}(C)\cup \mathrm{Int}(C')$ is a cycle. This new cycle contradicts the maximality of $C$.
\begin{flushright}
$\Box$
\end{flushright}

\begin{lem} \label{gug}
If a graph does not have any maximal cycle then this graph is a tree.
\end{lem}

\emph{Proof of the Lemma \ref{gug}}:

The proof is obvious. Indeed, if the graph does not have any cycle then there is nothing to prove, so we assume that the graph does have some cycles. In this case, the set of cycles is non-empty so it must have a maximal element for the partial order of definition \ref{defmax}. Thus the proof follows.
\begin{flushright}
$\Box$
\end{flushright}

\begin{lem} \label{eq}
Let $C$ be any cycle, then:
\begin{equation}
V_{out}^{C}-V_{in}^{C}=6.
\end{equation}
\end{lem}
\emph{Proof of Lemma \ref{eq}}:

This is a simple consequence of the Gauss-Bonnet Theorem. Let us travel along $C$, keeping the interior of $C$ on our left hand side. The cycle $C$
consist of finitely many line segments with exterior angles equal to $+60$ degrees at outgoing vertices and $-60$ degrees at incoming vertices. So
\begin{displaymath}
60V_{in}^{C}-60V_{out}^{C}=360.
\end{displaymath}
\begin{flushright}
$\Box$
\end{flushright}

\begin{cor} \label{impo}
Let $G$ be a critical graph which is not a tree. Then $G$ has at least $6$ boundary points.
\end{cor}

\emph{Proof of  Corollary \ref{impo}}:

By assumption, $G$ contains a maximal cycle $C$. Let $G'=G$ with the edges of $C$ removed. Then no two outgoing vertices of $C$ can be connected
in $G'$. Otherwise, let $v_1$, $v_2$ be outgoing vertices of $C$ connected in $G'$ by a path $\sigma$ of minimal length (over all paths and all pairs
$v_1$, $v_2$). Then $\sigma \cap C =\{v_1 ,v_2\}$. Let $\delta$ be one of the paths in $C$ which connects $v_1$ to $v_2$. Then $\sigma \cup \delta$ is a cycle, contradicting maximality of $C$.

Each connected component of $G'$ sitting outside $C$ is critical with boundary points consisting of a subset of boundary points of $G$ and exactly one outgoing vertex of $C$. It must have at least $2$ boundary points. Therefore
\begin{eqnarray*}
\vert \delta G\vert \geq V_{out}^{C}=6+V_{in}^{C}\geq 6.
\end{eqnarray*}

\begin{flushright}
$\Box$
\end{flushright}

\begin{lem} \label{forest}
Let $F$ be a disjoint union of $k$ $3$-boundary regular trees with a total of $n$ boundary points. Then 
\begin{eqnarray*}
\vert F\vert = 2n-2k.
\end{eqnarray*}
\end{lem}

\emph{Proof of Lemma \ref{forest}}:

For each component $T$ of $F$ with $n_{T}$ boundary points, $\vert T\vert=2n_{T}-2$. Summing over all components yields $\vert F\vert = 2n-2k$.
\begin{flushright}
$\Box$
\end{flushright}



\begin{lem} \label{kkk}
Let $G$ be a $3$-boundary regular critical graph with $n$ boundary points. If $G$ is a tree, then
\begin{eqnarray*}
\vert G\vert \leq 2n-2.
\end{eqnarray*}
\end{lem}
\emph{Proof of Lemma \ref{kkk}}:

The number of vertices of the binary tree with $n$ boundary points is equal to $2n-2$ ($n$ boundary points plus $n-2$ degree $3$ vertices). Therefore 
\begin{eqnarray*}
\vert G\vert \leq 2n-2.
\end{eqnarray*}
\begin{flushright}
$\Box$
\end{flushright}
\begin{lem} \label{onemaxcycle}
Let $G$ be a critical graph with $n$ boundary points. If $G$ has only one maximal cycle then $\vert G\vert \leq f_3(n-6)+2n$. Furtheremore, If $\vert G\vert=f_3(n-6)+2n$, then:
\begin{itemize}
\item Either $n=6$ and $G$ is the union of a $6$-cycle and its 6 outgoing edges,
\item or $n> 6$ and $G$ is obtained from a critical graph $G''$ attached to $n-6$ boundary points by completing its $n-6$ boundary points into a $2n-6$-cycle and adding the $n$ outgoing edges of this cycle.
\end{itemize}
\end{lem}
\emph{Proof of Lemma \ref{onemaxcycle}}:

$G$ has exactly $1$ maximal cycle $C$. Let $m'$ denote the number of outgoing vertices of the cycle. From Lemma \ref{eq}, the number of incoming vertices is equal to $m'-6$. Let $n'$ be the number of boundary points outside the cycle. Of course the number of boundary points inside the cycle is equal to $n-n'$.

If outside the cycle $C$, some (non maximal) cycle exists then the set of cycles outside $C$ must have a maximal element, hence a maximal cycle which is disjointed from the cycle $C$. This contradicts the assumption. Hence, outside $C$ the graph is a forest (a disjointed union of trees).

Now we remove all the edges of $C$, and consider vertices of $C$ as boundary points for the remaining graph $G'$. $G'$ consists of a graph
$G''$ whose edges were inside $C$ and of a collection $F$ of trees whose edges were outside $C$. $F$ has $m'+n'$ boundary points. 

Each outgoing vertex of $C$ is the root of one of the trees of the forest $F$. Thus the number of components of $F$ is $m'$. From Lemma \ref{forest}, $\vert F\vert= 2(m'+n')-2m'=2n'$. 

Each of these trees has at leat one boundary point, apart from its root, thus $n' \geq m'$.

$G''$ has at most $n-n'+m'-6$ boundary points. By the definition of the function $f_3(n)$ we know that the number of the vertices of $G''$ is at most $f_3(n+m'-n'-6)$. Then 
\begin{displaymath}
\vert G\vert \leq f_3(n+m'-n'-6)+2n'.
\end{displaymath}
Thus a $k \leq n-6$ exists such that
\begin{displaymath}
\vert G\vert \leq f_3(k)+2n'.
\end{displaymath}
On the other hand, we showed that $f_3(n)$ is non-decreasing. We conclude that 
\begin{equation} \label{se}
 \vert G\vert \leq f_3(n-6)+2n. 
\end{equation}
Equality implies that $n=n'=m'$ and that $F$ is a disjoint union of $n$ edges. Thus $F$ consists of the outgoing edges of $C$. If $G''$ is nonempty, the boundary points of $G''$ are the incoming vertices of $C$. Otherwise, $C$ has no incoming vertices, this means that $C$ is a $6$-cycle. This completes the proof of Lemma \ref{onemaxcycle}.
\begin{flushright}
$\Box$
\end{flushright}
\begin{lem} \label{moremaxcycles}
Let $G$ be a critical graph with $n$ boundary points. If $G$ has more than one maximal cycle, then $n'$ exists such that $ 6 \leq n' \leq n-4$ and
\begin{eqnarray*}
\vert G\vert \leq f_3(n')+f_3(n-n'+2)-2 .
\end{eqnarray*}
\end{lem}
\emph{Proof of Lemma \ref{moremaxcycles}}:

We shall use the following terminology:
\begin{de}
Let $C$ and $D$ be two maximal cycles in a critical $3$-boundary regular graph. A connecting set for $C$ and $D$ is a triple $(v_C, v_D ,\sigma)$ such that
$v_C$ (resp. $v_D$) is an outgoing vertex of $C$ (resp. $D$) and $\sigma$ a path which joins these vertices outside $\mathrm{C}\cup\mathrm{D}$.
\end{de}
Let $C$ and $D$ be two maximal cycles. Let $\sigma$ be a path connecting $C$ to $D$ with a minimum number of edges. Let $e$ be the first edge
traversed by $\sigma$. Then $e$ disconnects $C$ from $D$. Indeed, otherwise, there would exist a path connecting $C$ to $D$ away from $e$. Let $\gamma$ be the shortest path in $G\setminus\{e\}$ joining the endpoints of $e$. Then $\gamma \cup e$ is a cycle touching $C$ and thus contradicting maximality of $C$. So cutting the edge $e$ will disconnect $C$ from $D$.

Let $e'\subset e$ be a proper interval. Let $G'$(resp. $G''$) be the connected component of $G\setminus\{e\}$ containing $C$ (resp. $D$). Let $n'=\vert \delta G'\vert$. Then $G''$ has $n-n'+2$ boundary points (the cut through $e'$ produces two extra boundary points). By Lemma \ref{eq}, we conclude that $n' \geq 6$ and $n-n'+2 \geq 6$.

When computing the total number of vertices in $G$, the two extra boundary points created by cutting $e$ can be substracted off again. And so in final for an $n'$ such that $6 \leq n' \leq n-4$, we have
\begin{equation} \label{third}
\vert G\vert \leq f_3(n')+f_3(n-n'+2)-2.
\end{equation}
This completes the proof of Lemma \ref{moremaxcycles}.
\begin{flushright}
$\Box$
\end{flushright}

\begin{de}[Simple critical graphs]
Let $G$ be a critically embedded graph in the plane. Say that $G$ is {\em simple} if
\begin{itemize}
  \item Either $G$ is a tree with the following property : it does not contain paths consisting of 5 edges turning on the same side (like 5 consecutive edges of a convex hexagon).
\begin{center}
\includegraphics[width=2in]{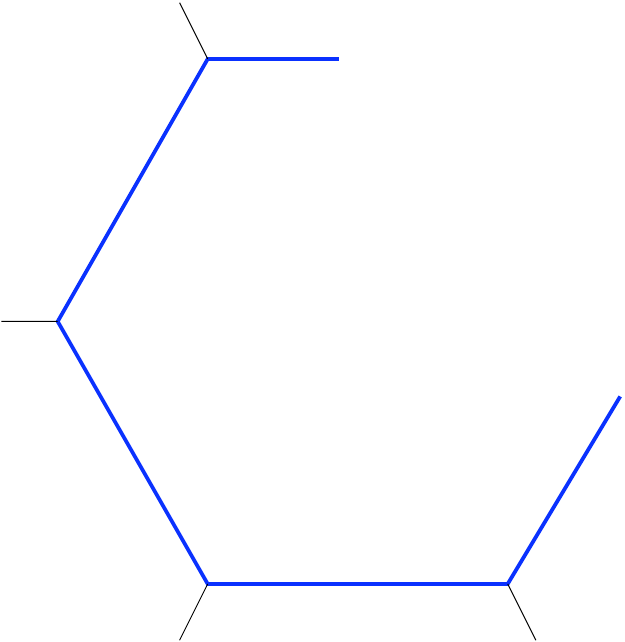}

Forbidden configuration in a simple critical tree.
\end{center}
  \item Or $G$ has a unique maximal cycle which surrounds all vertices except boundary points. Furthermore, no two consecutive vertices in the maximal cycle are both incoming vertices.
\begin{center}
\includegraphics[width=2in]{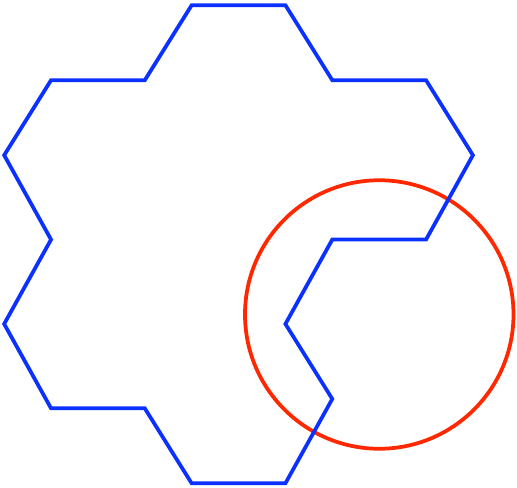}

Forbidden configuration in the maximal cycle of a simple critical graph
\end{center}
\end{itemize}
\end{de}
The pictures in section \ref{fig} all feature simple critical graphs. Note that on the set of boundary points of a simple critical graph, there is a natural anticlockwise order.

\begin{de}[Padding of a simple critical graph]
Let $G$ be a simple critically embedded graph in the plane. The {\em padding} of $G$, denoted by $P(G)$ is the critical graph obtained as follows.

We number $\delta G$ in anticlockwise order $a_1,\cdots,a_n$ where $a_{n+1}=a_1$. From each boundary point $a_i$ we draw two half-lines, $\alpha_{i}^{+}$ and $\alpha_{i}^{-}$ obtained by turning the edge which connects $a_i$ to $G$ by respectively $120$ and $240$ degrees. Next, one considers the portion of the cycle between $a_i$ and $a_{i+1}$, and completes this set of edges into an hexagon having two edges carried by $\alpha_{i}^{+}$ and $\alpha_{i+1}^{-}$. For this, one is led to place between 1 and 4 new vertices, depending on the configuration.
\begin{itemize}
\item If the angle $\angle (\alpha_{i}^{+},\alpha_{i+1}^{-})= \frac{-2\pi}{3}$, we cut the half-lines $\alpha_{i}^{+}$ and $\alpha_{i+1}^{-}$ on their intersection point $b_i$. We obtain two edges making an angle equal to $120$ degrees at $b_i$.
\begin{center}
\includegraphics[width=1in]{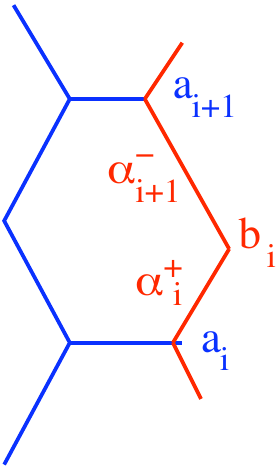}

Adding one vertex
\end{center}
\item  If the angle $\angle (\alpha_{i}^{+},\alpha_{i+1}^{-})= \frac{-\pi}{3}$, we place on $\alpha_{i}^{+}$ (resp $\alpha_{i+1}^{-}$), two points $b_{i}^{-}$ (resp $b_{i}^{+}$) such that the vector $\overrightarrow {b_{i}^{+}b_{i}^{-}}$ makes an angle equal to $\frac{\pi}{3}$ with $\alpha_{i}^{+}$ (i.e $\angle (\alpha_{i}^{+},b_{i}^{+}b_{i}^{-}) = \frac{\pi}{3}$ and $\angle (b_{i}^{+}b_{i}^{-}, \alpha_{i+1}^{-}) = \frac{\pi}{3}$).
\begin{center}
\includegraphics[width=2in]{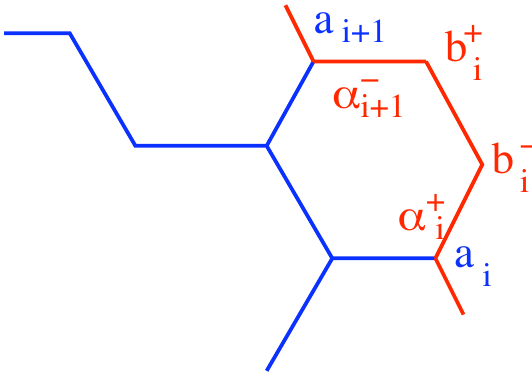}

Adding two vertices
\end{center}
\item If the angle $\angle (\alpha_{i}^{+},\alpha_{i+1}^{-})= 0$, we add a point $b_{i}^{-}$ on $\alpha_{i}^{+}$, a point $b_{i}^{+}$ on $\alpha_{i+1}^{-}$ and a point $_i$ such that $a_{i}b_{i}^{-}c_i b_{i}^{+}a_{i+1}\gamma_{i}$ is a hexagon with interior angles all equal to $120$ degree and where $\gamma_i$ is the vertex connected to  $a_{i}$ and $a_{i+1}$ by two edges of $G$.
\begin{center}
\includegraphics[width=1.5in]{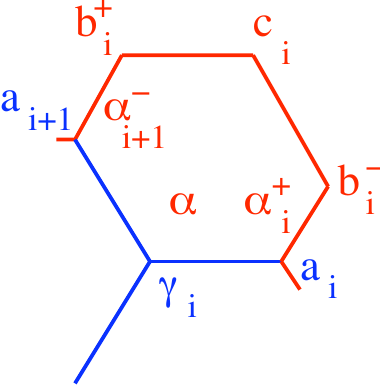}

Adding three vertices
\end{center}
\item If the angle $\angle (\alpha_{i}^{+},\alpha_{i+1}^{-})= \pi/3$ (this happens only if $G$ is the one-edge graph), we place a point $b_{i}^{-}$ on $\alpha_{i}^{+}$, a point $b_{i}^{+}$ on $\alpha_{i+1}^{-}$ and points $c_{i}^{-}$, $c_{i}^{+}$ such that $a_{i}b_{i}^{-}c_{i}^{-}c_{i}^{+}b_{i}^{+}a_{i+1}$ is a hexagon with interior angles all equal to $120$ degrees.
\begin{center}
\includegraphics[width=1.5in]{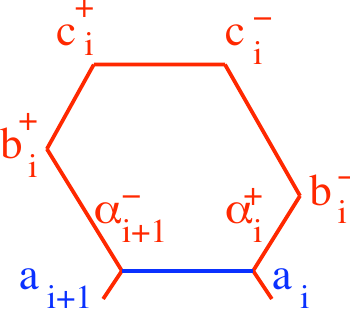}

Adding four vertices
\end{center}
\end{itemize}
The edges $\cdots \alpha_{i}^{+} (b_i^{+},\gamma_i,b_i^{-})\alpha_{i+1}^{-}\cdots $ form a cycle for which the vertices $a_i$ are incoming vertices and the $b_i$,$b_{i}^{+}$,$b_{i}^{-}$,$\gamma_i$ are outgoing vertices. We add to $G$ these edges with segments (and vertices) attached to each outgoing vertices  which will form the boundary points of the new critical graph $P(G)$.
\end{de}

The padded graph $P(G)$ is a simple critical graph. Indeed, by construction, it has a cycle which surrounds all vertices except boundary points. In this cycle, the $a_i$'s are incoming vertices, and between two consecutive $a_i$'s, outgoing vertices ($b_{i}^{\pm}$'s and $c_{i}^{\pm}$'s) are inserted. Therefore, the padding operation can be iterated. Of course one can give a much shorter definition for padding but we find our illustrative definition quite useful.



For every $n\geq 2$ we define a graph $H_n$ which is a subgraph of the standard tiling of the plane by regular hexagons.
\begin{de} [Graphs $H_n$]\label{defhn}
\begin{itemize}
\item $H_2$ consists of two vertices joined by a single edge.
\item For $2\leq n\leq 5$, $H_n$ consists  of a $3$-boundary regular tree with $n$ boundary points and $N(n)$ vertices.
\item $H_6$ is a critical graph consisting of a maximal cycle of length $6$ (a hexagon) and the outgoing edges and vertices attached to the hexagon.
\item $H_7$ is the critical graph which consists of adding a critical tree of length $3$ to a vertex of the hexagon of $H_6$.
\item For every $n\geq 8$, define $H_n$ inductively as the graph obtained by padding $H_{n-6}$, i.e $H_n=P(H_{n-6})$.
\end{itemize}
\end{de}

Section \ref{fig} shows the first $19$  $H_n$.

\emph{Remark} For every $n \geq 6$, $H_n$ has $n$ boundary points and only one maximal cycle of length $2n-6$. The $n$ boundary points correspond to the outgoing vertices of the cycle. 



\subsection{Proof of Theorem \ref{1}}
For every $n\geq 2$, let
\begin{eqnarray*}
N(n)=6k^2+6k+2l(k+1)-r,
\end{eqnarray*}
where $r=0$ if $n$ is a multiple of $6$ and $r=2$ otherwise ($N(n)$ is the expression given in Theorem \ref{1}).

We prove Theorem \ref{1} by showing that for every $n\geq 2$, $f_3(n)=N(n)$. First we list a few easy but important properties of $N(n)$.
\begin{lem} \label{prope}

\begin{itemize}
\item For every $n\geq 2$, $N(n)\geq 2n-2$ and for every $n\geq 6$, $N(n)>2n-2$.
\item For all $n\geq 8$, $N(n)=N(n-6)+2n$.
\item If $n\geq 2$ and $2\leq k\leq n$, $N(n)\geq N(k)+N(n-k+2)-2$ and more strongly if $n\geq 6$ and $3\leq k\leq n-1$, $N(n)>N(k)+N(n-k+2)-2$ except the cases $(n,k)=(7,3)$ or $(7,6)$.
\item For every $n\geq 2$, the number of vertices of $H_{n}$ equals $N(n)$.
\end{itemize}
\end{lem}
Proof of the first two parts are straightforward. For the third part one can use the following easily-verified inequalities: 
\begin{eqnarray*}
\frac{1}{6}n^2+n-\frac{7}{6}\leq N(n)\leq \frac{1}{6}n^2+n.
\end{eqnarray*}

For the final part, the proof is settled using the recursive definition of $H_n$.
\begin{flushright}
$\Box$
\end{flushright}

We prove Theorem \ref{1} by induction. We show that for every $n\geq 2$, $f_3(n)=N(n)$.

\emph{Basis of the Induction}:

Corollary \ref{impo} and Lemma \ref{gug} imply that for $2\leq n\leq 5$ the extremal graph is a tree. It is then straightforward that the maximum number of vertices of such a tree equals $N(n)$. Morever, the isomorphism classes of such trees are unique. For $n=6$ the required equality follows from the Lemma \ref{onemaxcycle} and for $n=7$ the next Lemma provides the desired equality.


\begin{lem} \label{last}
The maximum number of vertices of a $3$-boundary regular critical graph with $7$ boundary points is equal to $N(7)=14$. Furthermore, every $3$-boundary regular critical graph with $7$ boundary points and $14$ vertices is isomorphic to $H_7$.
\end{lem}
\emph{Proof of Lemma \ref{last}} 

$H_7$ can be presented in the tiled plane by a hexagon with $5$ segments attached to $5$ vertices of the hexagon and a tree with $4$ vertices having one boundary point in the $6$th vertex of the hexagon.

We need to prove that this graph has the maximum number of vertices among all $3$-boundary regular critical graphs with $7$ boundary points.

Let's suppose that a graph $H$ exists having more vertices than $H_7$. According to Corollary \ref{impo}, $H$ has only one maximal cycle and the length of the maximal cycle is equal to $6$. Then each vertex of the maximal cycle can be considered as a boundary point for a critical tree (attached to the vertex). If there exist two vertices of the maximal cycle such that the two critical trees attached to them are not segments, then the number of boundary points of $H$ will be more than $7$ and this is not possible. So to $5$ vertices of the maximal cycle are attached to $5$ segments, and the proof of the Lemma follows.
\begin{flushright}
$\Box$
\end{flushright}

The above discussion together with the Lemma \ref{last} shows that for $2\leq n\leq 7$ we have $f_3(n)=N(n)$ and the unique extremal example is $H_n$.


We now prove Theorem \ref{1} for $n\geq 8$ with the inductive assumption including the fact that the unique extremal example showing $f_3(n')=N(n')$ is $H_{n'}$ for $n'<n$. The graph maximising the number of vertices with $n$ boundary points is either a tree or has one and only one maximal cycle or has more than one maximal cycle.
From the reasoning of Lemma \ref{onemaxcycle} and Lemma \ref{moremaxcycles} together with the inequalities involving $N(n)$ in the Lemma \ref{prope} it follows that for $n\geq 8$ we have $f_3(n)\leq N(n)$ and if $f_3(n)=N(n)$ then the extremal example must arise when there is one maximal cycle by padding the extremal example for $n-6$, so it must be $H_n$. This finalises the Proof of Theorem \ref{1}. 
\begin{flushright}
$\Box$
\end{flushright}

\section{$4$-Boundary Regular Critical Graphs}

Here we prove Theorem \ref{sec}.

Let $G$ be a 4-boundary regular critical graph with $n$ boundary points. Then $G$ is made up of line segments intersecting each other in a way that when any two segments intersect at a point (vertex of $G$), there are no other segments passing through the intersecting point. Thus every intersection point will be a vertex of $G$ and the critical condition is verified. Every line segment joins two of the boundary points. Then the problem of estimating $f_4(n)$ is equivalent to finding the maximum number of intersection points of $n/2$ line segments in the plane such that only two lines pass through the intersecting points.

For $n$ odd it is impossible to attach a critical graph of degree $4$ to $n$ points, and $f_4(n)=0$ in this case. For $n$ even, the number of intersecting points will not exceed $\binom{n/2}{2}$.

To complete the proof of the theorem, we show that for every even $n$, there exists a collection of $n/2$ line segments intersecting at exactly $\binom{n/2}{2}$ points. We prove the existence of such a collection by induction on $n$. For $n=1,2$ it is obvious. We suppose that such a collection is constructed for $n$ and we need to add a single line $L$ to this collection such that $L$ does not pass through the intersection points, and such that $L$ intersects all the lines of the collection. As the number of lines and their intersections is finite, it is always possible to add such a line $L$ with the required property. From the existence of such a collection the proof of Theorem \ref{sec} follows.
\begin{flushright}
$\Box$
\end{flushright}

\section{Remarks and Open Questions}
The problem of estimating the maximum number of vertices of a critical graph attached to some boundary points in the plane for the case of $3$ and $4$-boundary regular graphs turned out to be very elementary. We saw that without any difficulties, we could even classify maximising graphs. But the same question for a non-necessarily boundary regular graph is more complicated. Indeed, for vertices of degree $3$ the angle between any two outgoing edges is equal to $120$ degrees. This simple fact allows us to have a Gauss-Bonnet type lemma, and makes our estimates possible. But for degrees greater than $4$, there are infinite possible geometric configurations of outgoing edges.

\begin{nott}
Let $d$ be a natural number, we denote by $g_{d}(n)$ the supremum of the number of vertices of critical graphs with $n$ boundary points and degree bounded by $d$.
\end{nott}

The general questions are
\begin{itemize}
\item for $d=4$, find a better upper bound for $g_{4}(n)$.
\item for $d>4$, is $g_{d}(n)$ finite ?
\end{itemize}
An initial guess would be that $g_{d}(n)=f_3(n)$ and that the class of $3$-boundary regular critical graphs have the largest $g_{d}(n)$ for all the value of $n$ and among all the critical graphs with bounded degrees. Indeed one can imagine that locally, every critical graph having a vertex of degree greater than $3$  can be mapped by a homotopy to a $3$-boundary regular critical graph, such that the number of vertices of the image by the homotopy increase. The non-obvious part is that these local homotopies will move the position of the boundary points and that we can't re-assemble the local part of the graph correctly together and get a new $3$-boundary regular critical graph. Thus the initial guess may be misleading.

In the introduction, we briefly discussed graphs which minimise the length functional. It is interesting to note that all of the edges of the $3$-boundary regular graphs maximising the number of vertices with $n$ boundary points have length equal to $1$. The question which is not clear is the following : Is it true that every edge of a general critical graph with $n$ boundary points maximising the number of vertices has a length equal to $1$? Of course, it's easy to construct critical graphs with a finite number of boundary points with different lengths of edges- but what about the ones maximising the number of vertices? 

In order to study the variational Plateau problem mentioned in the introduction, it is important to define a topology in the space of graphs with fixed boundary points. This variational problem is best formulated in the language of geometric measure theory, where for instance, one defines the space of currents or varifolds. Every planar graph can be considered as a $1$-current or a $1$-varifold. One could study different interesting functionals on such spaces and seek out graphs minimising those functionals. Allard and Almgren gave an example of a family of $3$-boundary regular (weighted) minimal graphs with $16$ boundary points and with arbitrarily large numbers of vertices. The minimality condition is with regards to a functional different from the one defined in the introduction. These examples are known as the \emph{spiderweb-like varifolds}, see \cite{Alm1}, \cite{Alm2} and \cite{Alm3}. They motivate the following conjecture:

\begin{conj}
There exist critical graphs in the plane with a finite number of boundary points and arbitrarily large number of vertices.
\end{conj}

If this conjecture is true, it will be interesting to study infinite critical graphs. This could also motivate the study of Morse theory in the
space of $1$-cycles with infinitely many edges.

However, paradoxically:
\begin{conj}
For all $d \geq 3$, the number of vertices of a $d$-boundary regular critical graph with $n$ boundary points is $ \leq f_3 (n)$.
\end{conj}
One can begin with the case in which all of the angles of the $d$-boundary regular critical graphs are equal to $\frac{2\pi}{d}$ and try to obtain an intermediate result like Lemma \ref{eq}.

The problem of estimating the maximum number of vertices in a critical graph with some fixed conditions can also be asked in a more general context. One may also ask the same questions about the graphs embedded in Riemannian manifolds. We shall generalise the definition of critical planar graphs to :

\begin{de} [Critical graphs in Riemannian manifolds]
Let $M$ be a Riemannian manifold and let $A$ be a finite subset of $M$. A \emph{critical} graph with boundary points $A$ is a finite embedded graph $G$ in $M$ such that the following conditions are satisfied:
\begin{enumerate}
\item Each edge of the graph is a geodesic segment.
\item Every $a \in A$ represents a vertex of degree 1 in $G$.
\item The sum of unit vectors of edges outgoing  from each vertex of degree greater than 1 is equal to zero.
\end{enumerate}
\end{de}

The known results concern mostly the $2$-sphere with a non necessarily canonical metric (see \cite{mo1} and \cite{he}). In \cite{he}, the author studies the critical graphs on the $2$-sphere and shows that the $2$-sphere contains critical graphs with arbitrarily large number of vertices. The same question is interesting if one imposes that all the vertices and boundary points have to lie on a hemi-sphere (this is equivalent to studying the same problem for the plane with a different Riemannian metric). The case of a flat cylinder is also interesting, since a flat cylinder contains arbitrarily large $d$-boundary regular graphs with at most $6$ attaching points. It is interesting to give a value for $g_d(n)$ in the Hyperbolic plane and compare the result to the case of plane or (hemi-)sphere (at least for $d=3,4$).

\newpage
\section{Figures}
\label{fig}

\begin{center}\includegraphics[width=1.5in]{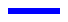}
\hspace{-1cm}\includegraphics[width=1.5in]{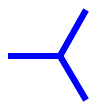}
\hspace{-1cm}\includegraphics[width=1.5in]{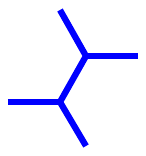}
\hspace{-1cm}\includegraphics[width=1.5in]{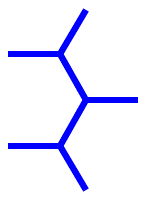}
\hspace{-1cm}\includegraphics[width=1.5in]{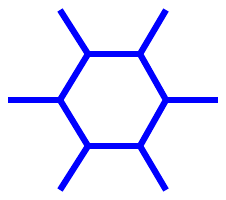}
\end{center}
\vspace{-2cm}
\begin{center}
\includegraphics[width=1.5in]{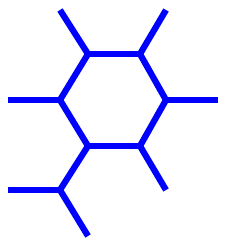}
\hspace{-1cm}\includegraphics[width=1.5in]{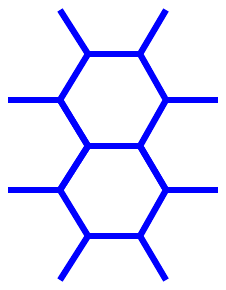}
\hspace{-1cm}\includegraphics[width=1.5in]{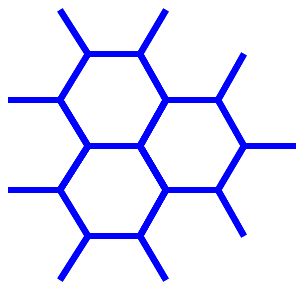}
\hspace{-1cm}\includegraphics[width=1.5in]{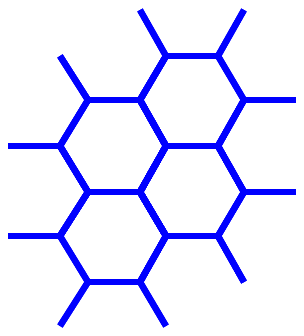}
\hspace{-1cm}\includegraphics[width=1.5in]{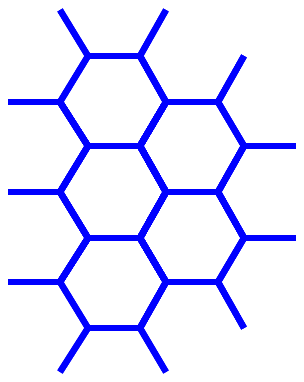}
\end{center}
\vspace{-2cm}
\begin{center}
\includegraphics[width=1.5in]{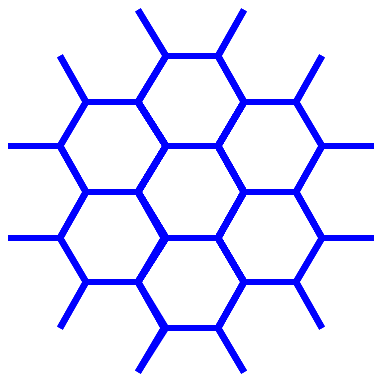}
\hspace{-1cm}\includegraphics[width=1.5in]{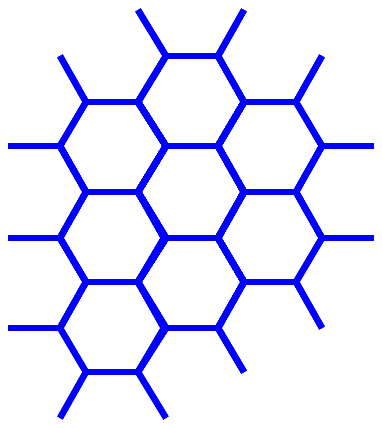}
\hspace{-1cm}\includegraphics[width=1.5in]{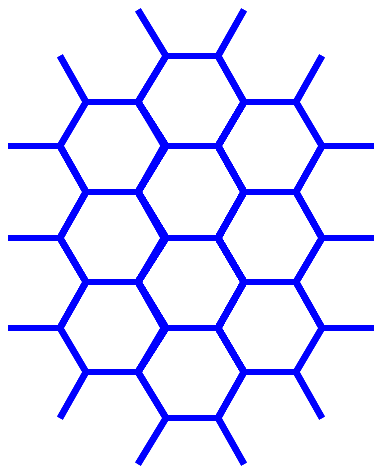}
\hspace{-1cm}\includegraphics[width=1.5in]{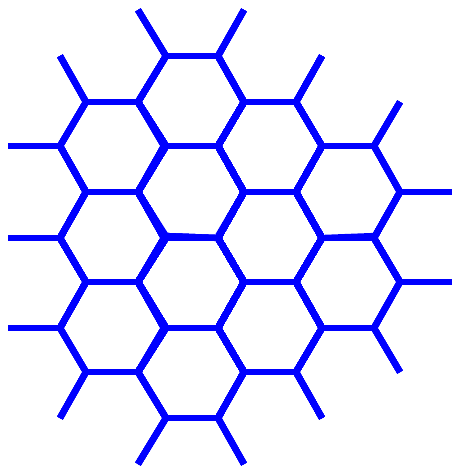}
\hspace{-1cm}\includegraphics[width=1.5in]{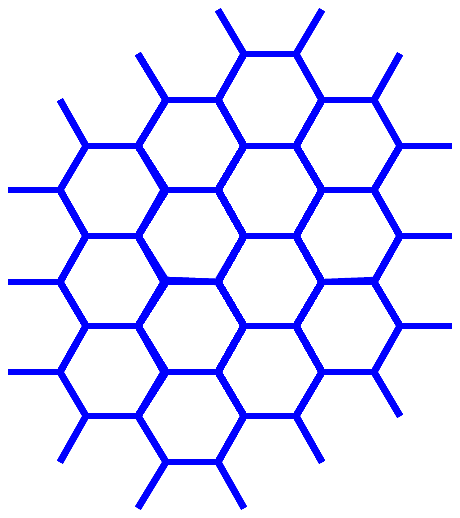}

\end{center}
\vspace{-2cm}
\begin{center}
\includegraphics[width=1.5in]{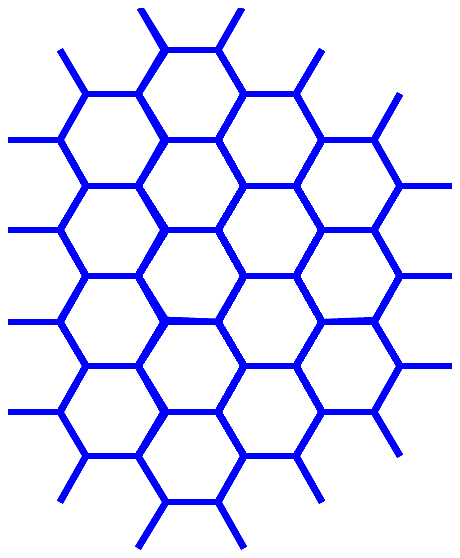}
\hspace{-1cm}\includegraphics[width=1.5in]{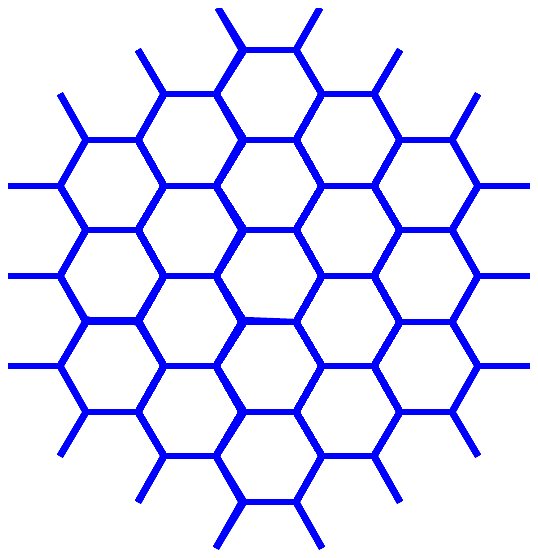}
\hspace{-1cm}\includegraphics[width=1.5in]{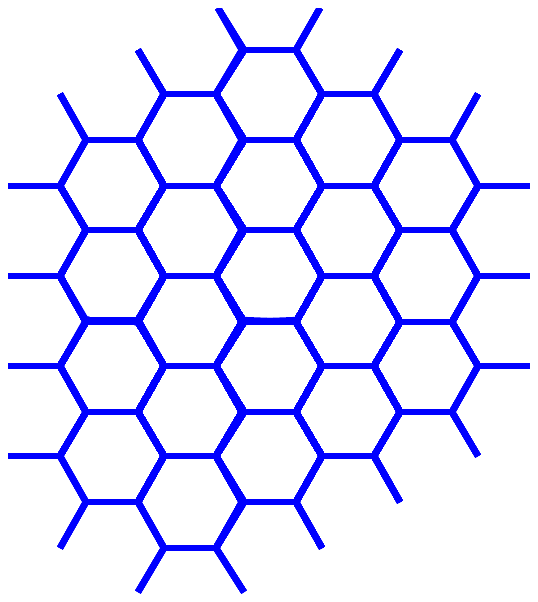}
\hspace{-1cm}\includegraphics[width=1.5in]{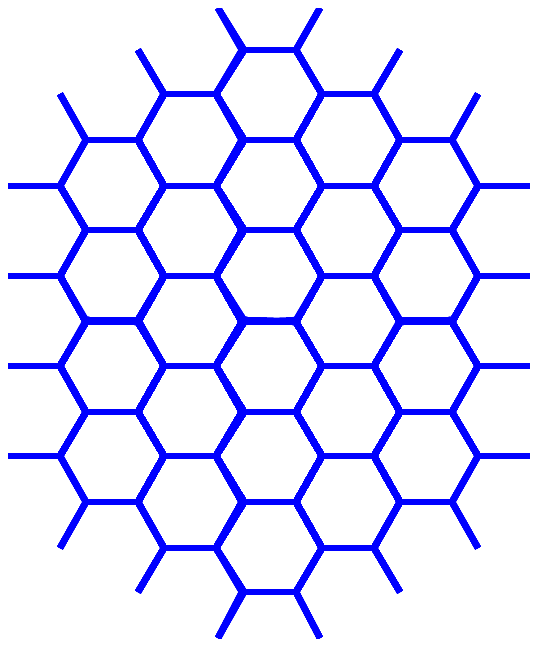}
\end{center}

\section{Acknowledgments}

I am grateful to M.Gromov for suggesting this interesting problem to me in connection with his recent
paper \cite{gromov}, and to P.Pansu for helping me realise this work.

\bibliographystyle{plain}
\bibliography{mabiblio}
\end{document}